\documentclass[reprint,amsmath,amssymb,pre]{revtex4-1}

\usepackage{amsfonts}
\usepackage{amsthm}
\usepackage{graphicx}
\usepackage{xspace}
\usepackage{tikz}
\usetikzlibrary{calc}
\usetikzlibrary{matrix}
\usetikzlibrary{patterns}
\usetikzlibrary{arrows}

\newtheorem{theorem}{Theorem}

\newtheorem{definition}{Definition}

\providecommand{\comment}[1] {}

\newcommand*{\sect}[1]	{$\S\,#1$}
\newcommand*{\latin}[1]	{\emph{#1}}
\newcommand*{\ie}			{\latin{i.e.}\@\xspace}
\newcommand*{\eg}			{\latin{e.g.}\@\xspace}
\newcommand*{\cf}			{\latin{cf.}\@\xspace}

\newcommand*{\etal}		{\latin{et al.}\@\xspace}
\newcommand*{\figrf}[1]	{fig.~\ref{f:#1}}
\newcommand*{\Figrf}[1]	{Fig.~\ref{f:#1}}
\newcommand*{\eqrf}[1]	{(\ref{eq:#1})}
\newcommand*{\secrf}[1]	{\sect{\ref{s:#1}}}

\providecommand{\figplacement}{tbp}
\newenvironment{fig}[2][\figplacement]
	{%
		\newcommand{\bbfig}[1]{\begin{figure}[##1]}
		\expandafter\bbfig\expandafter{#1}
		\def\fglbl{f:#2}
		\centering
	}{%
		\label{\fglbl}
		\end{figure}
	}

\newcommand*{\define}	{\overset{\text{\tiny def}}{=}}
\newcommand*{\integers} 	{\mathbb{Z}}

\newcommand*{\card}[1]		{\abs{#1}}

\newcommand*\half				{\tfrac{1}{2}}
\newcommand*{\tends}			{\rightarrow}
\renewcommand*{\vec}[1]		{\mathbf{#1}}
\providecommand{\abs}[1]	{\lvert#1\rvert}

\newcommand\X{\mathcal{X}}

\newcommand\A{\mathcal{A}}

\newcommand\sps{\,.\,}
\newcommand\Ipred{\mathcal{I}_{\mathrm{pred}}}

\newcommand\IXZ{\overline{\underline{\mathcal{I}}}}

\renewcommand\define{\triangleq}

\newcommand\cmin{\!-\!}
\newcommand\modulo[2]{#1\!\!\!\!\!\mod#2}

\newcommand{\past}[1]{\loarrow{#1}}
\newcommand{\fut}[1]{\roarrow{#1}}
\newcommand\parity[2]{P^{#1}_{2,#2}}

\newcommand\iftwocol[2]{#1}
\begin{document}
	\title{A measure of statistical complexity based on predictive information}
	\date{\today}
	\author{Samer A. Abdallah}
	\author{Mark D. Plumbley}
	\affiliation{Queen Mary University of London}

	\begin{abstract}
	We introduce an information theoretic measure
	of statistical structure, called `binding information', for sets of 
	random variables, and compare it
	with several previously proposed measures
	including excess entropy,
	Bialek \etal's predictive information,
	and the multi-information.
	We derive some of the properties of the binding 
	information, 
	particularly in relation to the multi-information,
	and show that, 
	for finite sets of
	binary random variables, 
	the processes which maximises
	binding information are the `parity' processes.
	Finally we discuss some of the implications this has for the
	use of the binding information as a measure of complexity.
	\end{abstract}
	\maketitle

	\section{Introduction}
	\label{s:intro}

	The concepts of `structure', `pattern' and `complexity' are relevant in many
	fields of inquiry: physics, biology, cognitive sciences, machine learning, the arts
	and so on; but are vague enough to resist being quantified in a single definitive
	manner. One approach, which we adopt here, is to attempt to characterise them in statistical
	terms, for \emph{distributions} over configurations of some system,
	using the tools of information theory \cite{CoverThomas}. 

	In this letter, we propose a measure of statistical
	structure
	based on the concept of
	\emph{predictive information rate} (PIR) \cite{AbdallahPlumbley2009},
	which measures an aspect of temporal
	dependency not captured by previously proposed measures. We 
	review a number of these earlier proposals and the PIR, and then
	define the \emph{binding information} as the extensive counterpart of the PIR 
	applicable to arbitrary countable sets of random variables.
	After describing some of its properties, we identify some finite
	discrete processes that maximise the binding information.

	In the following,
	if $X$  is a random process
	indexed by a set $\A$, and $\mathcal{B}\subseteq\A$,
	then $X_{\mathcal{B}}$ denotes the compound random variable
	(random `vector') formed by taking $X_\alpha$ for each $\alpha \in \mathcal{B}$. 
	The set of integers from $M$ to $N$ inclusive will be written $M..N,$
	and $\setminus$ will denote the set difference operator, 
	so, for example, $X_{1..3\setminus\{2\}} \equiv (X_1,X_3)$.

\section{Background}
	Suppose that  $(\ldots,X_{-1},X_0,X_1,\ldots)$ is a bi-infinite stationary sequence of
	random variables,
	and that
	$\forall t\in\integers$, 
	the random variable $X_t$ takes values in a discrete set $\X$. Let $\mu$ be
	the associated shift-invariant probability measure. Stationarity implies
	that the probability distribution associated with any
	contiguous block of $N$ variables $(X_{t+1},\ldots,X_{t+N})$ 
	is independent of $t$, and therefore 
	we can define a shift-invariant block entropy function: 
	\begin{equation}
		\label{eq:block-entropy}
		H(N) \define H(X_1,\ldots,X_N) = 
			 \!\!\sum_{\vec{x}\in\X^N}\!\! -p_\mu^N(\vec{x}) \log p_\mu^N(\vec{x}),
	\end{equation}
	where $p_\mu^N : \X^N \to [0,1]$ is the unique probability mass function
	for any $N$ consecutive variables in the sequence,
	$p_\mu^N(\vec{x}) \define \Pr(X_1=x_1\land\ldots\land X_N=x_N)$.

\comment{
	\begin{fig}{venn-example}
		\newcommand\rad{2em}%
		\newcommand\circo{circle (3.1em)}%
		\newcommand\labrad{4.1em}
		\newcommand\bound{(-6em,-5em) rectangle (6em,6em)}
		\newcommand\colsep{\ }
		\newcommand\clipin[1]{\clip (#1) \circo;}%
		\newcommand\clipout[1]{\clip \bound (#1) \circo;}%
		\newcommand\cliptwo[3]{%
			\begin{scope}
				\clipin{#1};
				\clipin{#2};
				\clipout{#3};
				\fill[black!30] \bound;
			\end{scope}
		}%
		\newcommand\clipone[3]{%
			\begin{scope}
				\clipin{#1};
				\clipout{#2};
				\clipout{#3};
				\fill[black!15] \bound;
			\end{scope}
		}%
		\begin{tabular}{c@{\colsep}c}
			\begin{tikzpicture}[baseline=0pt]
				\coordinate (p1) at (90:\rad);
				\coordinate (p2) at (210:\rad);
				\coordinate (p3) at (-30:\rad);
				\clipone{p1}{p2}{p3};
				\clipone{p2}{p3}{p1};
				\clipone{p3}{p1}{p2};
				\cliptwo{p1}{p2}{p3};
				\cliptwo{p2}{p3}{p1};
				\cliptwo{p3}{p1}{p2};
            \begin{scope}
               \clip (p1) \circo;
               \clip (p2) \circo;
               \clip (p3) \circo;
               \fill[black!45] \bound;
            \end{scope}
				\draw (p1) \circo;
				\draw (p2) \circo;
				\draw (p3) \circo;
				\path 
					(barycentric cs:p3=1,p1=-0.2,p2=-0.1) +(0ex,0) node {$I_{3|12}$}
					(barycentric cs:p1=1,p2=-0.2,p3=-0.1) +(0ex,0) node {$I_{1|23}$}
					(barycentric cs:p2=1,p3=-0.2,p1=-0.1) +(0ex,0) node {$I_{2|13}$}
					(barycentric cs:p3=1,p2=1,p1=-0.55) +(0ex,0) node {$I_{23|1}$}
					(barycentric cs:p1=1,p3=1,p2=-0.55) +(0ex,0) node {$I_{13|2}$}
					(barycentric cs:p2=1,p1=1,p3=-0.55) +(0ex,0) node {$I_{12|3}$}
					(barycentric cs:p3=1,p2=1,p1=1) node {$I_{123}$}
					;
				\path
					(p1) +(140:\labrad) node {$X_1$}
					(p2) +(-140:\labrad) node {$X_2$}
					(p3) +(-40:\labrad) node {$X_3$};
			\end{tikzpicture}
			&
			\parbox{0.5\linewidth}{
				\small%
				\begin{align*}
					I_{1|23} &= H(X_1|X_2,X_3) \\
					I_{13|2} &= I(X_1;X_3|X_2) \\
					I_{1|23} + I_{13|2} &= H(X_1|X_2) \\
					I_{12|3} + I_{123} &= I(X_1;X_2) 
				\end{align*}
			}
		\end{tabular}
		\caption{
		Venn diagram visualisation of entropies and mutual informations
		for three random variables $X_1$, $X_2$ and $X_3$. The areas of each
		of the circles represent $H(X_1)$, $H(X_2)$ and $H(X_3)$ respectively.
		The total shaded area is the joint entropy $H(X_1,X_2,X_3)$.
		The central area $I_{123}$ is the co-information.
		Some other information measures are indicated in the legend.
		}
	\end{fig}
}

	The \emph{entropy rate} $h_\mu$ 
	has two
	equivalent definitions in terms of the block entropy function \cite[Ch. 4]{CoverThomas}:
	\begin{equation}
		\label{eq:entro-rate-alt}
		h_\mu \define \lim_{N\tends\infty} \frac{H(N)}{N}
			= \lim_{N\tends\infty} H(N) - H(N-1).
	\end{equation}
	The block entropy function can also be used to express
	the mutual information between two contiguous segments of the
	sequence 
	of length $N$ and $M$ respectively:
	\begin{equation}
		I(X_{-N..-1};X_{0..M-1}) = H(N) + H(M) - H(N+M).
	\end{equation}
	If we let both block lengths $N$ and $M$ tend to infinity,
	we obtain what has been called the \emph{excess entropy}
	\cite{CrutchfieldPackard1983} or the 
	\emph{effective measure complexity} \cite{Grassberger1986}.
	It is the
	amount of information about the infinite future that can be obtained, 
	on average, by observing the infinite past:
	\begin{equation}
		E = \lim_{N\tends\infty} 2H(N)-H(2N).
	\end{equation}
\comment{
	It can also be expressed in terms of the $h_\mu(N)$ defined by 
	Crutchfield and Feldman \cite{CrutchfieldFeldman1997} as 
	$h_\mu(N) \define H(X_N|X_{1:N-1}) = H(N)-H(N-1)$, which can be thought of
	as an estimate of the entropy rate obtained from the finite
	dimensional marginal distribution $p_\mu^N$. 
	In terms of $h_\mu(\cdot)$, the excess entropy can be
	expressed as \cite{CrutchfieldFeldman1997}
	\begin{equation}
		\label{eq:excess-entropy}
		E = \sum_{M=1}^\infty (h_\mu(M) - h_\mu).
	\end{equation}
	Grassberger \cite{Grassberger1986} and others
	\cite{LindgrenNordahl1988,Li1991}
	have commented on the manner in 
	which $h_\mu(N)$ approaches its
	limit $h_\mu$, noting that in certain types of random process with
	long-range correlations, the convergence can be so slow
	(as a power of $N$ rather than exponential) that the excess entropy
	is infinite, and that this is indicative of a certain kind of
	complexity. This phenomenon was examined in more detail
	by Bialek \etal \cite{BialekNemenmanTishby2001}, who defined
}
	Bialek \etal \cite{BialekNemenmanTishby2001} defined
	the \emph{predictive information} $\Ipred(N)$ as
	the mutual information between a 
	block of length $N$ and the
	infinite future following it:
	\begin{equation}
		\label{eq:predinfo}
		\Ipred(N) \define \lim_{M\tends\infty} H(N) + H(M) - H(N+M).
	\end{equation}
	They
	showed that even if 
	$\Ipred(N)$ diverges as $N$ tends to infinity, the \emph{manner}
	of its divergence reveals something about 
	the learnability of the underlying random process.
	Bialek \etal 
	also emphasised that $\Ipred(N)$ is the \emph{sub-extensive}
	component of the entropy: 
	if $N h_\mu$ is the purely extensive (\ie, linear in $N$) component of the 
	entropy, then $\Ipred(N)$ is the difference between
	the block entropy $H(N)$ and its extensive component:
	\begin{equation}
		\label{eq:subextensive-entropy}
		H(N) = N h_\mu + \Ipred(N).
	\end{equation}

	The \emph{multi-information} \cite{StudenyVejnarova1998} is 
	defined for any collection of $N$ random variables $(X_1,\dots,X_N)$ as
	\begin{equation}
		\label{eq:multi-info}
		I(X_{1..N}) \define - H(X_{1..N}) + \sum_{i\in 1..N} H(X_i). 
	\end{equation}
	For $N=2$, the multi-information reduces to the mutual information
	$I(X_1;X_2)$, while for $N>2$, $I(X_{1:N})$ continues to be a measure
	of dependence, being zero if and only if the variables
	are statistically independent. 
	In the thermodynamic limit, the intensive \emph{multi-information rate}
	(\cf Dubnov's \emph{information rate} \cite{Dubnov2004}) can be defined as
	\begin{equation}
		\rho_\mu \define \lim_{N\tends\infty} I(X_{1..N}) - I(X_{1..N-1}).
	\end{equation}
	It can easily be shown that
	$\rho_\mu = \Ipred(1) = H(1) - h_\mu$.
	Erb and Ay \cite{ErbAy2004} studied this quantity (they
	call it $I$) and 
	showed that, in the present terminology,
	\begin{equation}
		\label{eq:subextensive-multiinfo}
		I(X_{1..N}) + \Ipred(N) = N\rho_\mu.
	\end{equation}
	Comparing this with \eqrf{subextensive-entropy}, we see that 
	$\Ipred(N)$ is also the sub-extensive component of the multi-information.
	Thus, all of the measures considered so far,
	being linearly dependent in various ways,
	are closely related.

	Another class of measures,
	including Grassberger's \emph{true measure complexity}
	\cite{Grassberger1986}
	and Crutchfield \etal's \emph{statistical complexity} $C_\mu$
	\cite{CrutchfieldYoung1989,CrutchfieldFeldman1997},
	is based on the properties
	of stochastic automata that model the process under consideration.
	These have some interesting properties 
	but are beyond the scope of this letter.

 	\begin{figure}
		\centering
		\newcommand\rad{1.6em}%
		\newcommand\ovoid[1]{%
			++(-#1,\rad) 
			-- ++(2 * #1,0em) arc (90:-90:\rad)
 			-- ++(-2 * #1,0em) arc (270:90:\rad) 
		}%
		\newcommand\axis{2.75em}%
		\iftwocol{%
			\newcommand\olap{0.85em}%
			\newcommand\offs{3.6em}
			\newcommand\colsep{\hspace{5em}}
		}{%
			\newcommand\olap{1em}%
			\newcommand\offs{3.75em}
			\newcommand\colsep{\qquad}
		}%
		\newcommand\longblob{\ovoid{\axis}}
		\newcommand\shortblob{\ovoid{1.75em}}
		\begin{tabular}{c}
			\newcommand\thepast{\shortblob}
			\newcommand\future{\longblob}
         \begin{tikzpicture}
            \coordinate (p3) at (-\offs,0em);
            \coordinate (p1) at (-\offs+1em,0em);
            \coordinate (p2) at (\offs,0em);
            \begin{scope}
               \clip (p1) \thepast;
               \clip (p2) \future;
               \fill[lightgray] (-1,-1) rectangle (1,1);
            \end{scope}
            \draw[white] (p3) \longblob;
            \draw (p1) +(-0.75em,0em) node{\shortstack{finite\\past}} \thepast;
            \draw (p2) +(0.5em,0em) node{\shortstack{infinite\\future}} \future;
            \path (0,0) node (future) {\small PI};
            \path (p1) +(-0.5em,\rad) node [anchor=south] {$X_{-N:-1}$};
            \path (p2) +(2em,\rad) node [anchor=south] {$X_0,\ldots$};
         \end{tikzpicture} \\
			(a) predictive information $\Ipred(N)$ 
         \\[0.75em]
			\newcommand\blob{\longblob}
         \begin{tikzpicture}
            \coordinate (p1) at (-\offs,0em);
            \coordinate (p2) at (\offs,0em);
            \begin{scope}
               \clip (p1) \blob;
               \clip (p2) \blob;
               \fill[lightgray] (-1,-1) rectangle (1,1);
            \end{scope}
            \draw (p1) +(-0.5em,0em) node{\shortstack{infinite\\past}} \blob;
            \draw (p2) +(0.5em,0em) node{\shortstack{infinite\\future}} \blob;
            \path (0,0) node (future) {$E$};
            \path (p1) +(-2em,\rad) node [anchor=south] {$\ldots,X_{-1}$};
            \path (p2) +(2em,\rad) node [anchor=south] {$X_0,\ldots$};
         \end{tikzpicture} \\
         (b) excess entropy $E$ 
			\\[0.75em]
\comment{
         \begin{tikzpicture}
            \coordinate (p1) at (-\offs,0em);
            \coordinate (p2) at (\offs,0em);
				\path[fill=lightgray,draw=black] (-\olap,0) circle (\rad);
            \draw (p1) +(-1.75em,0em) node{\shortstack{infinite\\past}} \longblob;
            \draw (p2) +(0.5em,0em) node{\shortstack{infinite\\future}} \longblob;
				\path (-1.8em,0) node {$C_\mu$};
            \path (p1) +(-2em,\rad) node [anchor=south] {$\ldots,X_{-1}$};
            \path (p2) +(2em,\rad) node [anchor=south] {$X_0,\ldots$};
         \end{tikzpicture}
			 \\ (c) complexity $C_\mu$ 
         \\
}
         \begin{tikzpicture}[baseline=-1em]
				\newcommand\rc{2em}
				\newcommand\throw{2.5em}
            \coordinate (p1) at (210:1.5em);
            \coordinate (p2) at (90:0.7em);
            \coordinate (p3) at (-30:1.5em);
				\newcommand\bound{(-7em,-3.0em) rectangle (7em,3.5em)}
				\newcommand\present{(p2) circle (\rc)}
				\newcommand\thepast{(p1) ++(-\throw,0) \ovoid{\throw}}
				\newcommand\future{(p3) ++(\throw,0) \ovoid{\throw}}
				\newcommand\fillclipped[2]{%
					\begin{scope}[even odd rule]
						\foreach \thing in {#2} {\clip \thing;}
						\fill[black!#1] \bound;
					\end{scope}%
				}%
				\fillclipped{30}{\present,\future,\bound \thepast}
				\fillclipped{15}{\present,\bound \future,\bound \thepast}
            \draw \future;
            \fillclipped{45}{\present,\thepast}
            \draw \thepast;
            \draw \present;
            \node at (barycentric cs:p2=1,p1=-0.17,p3=-0.17) {$r_\mu$};
            \node at (barycentric cs:p1=-0.4,p2=1.0,p3=1) {$b_\mu$};
            \node at (barycentric cs:p3=0,p2=1,p1=1) [shape=rectangle,fill=black!45,inner sep=1pt]{$\rho_\mu$};
				\path (p2) +(140:3em) node {$X_0$};
            \path (p3) +(3em,0em) node  {\shortstack{infinite\\future}};
            \path (p1) +(-3em,0em) node  {\shortstack{infinite\\past}};
            \path (p1) +(-4em,\rad) node [anchor=south] {$\ldots,X_{-1}$};
            \path (p3) +(4em,\rad) node [anchor=south] {$X_1,\ldots$};
         \end{tikzpicture} \\
			(c) predictive information rate $b_\mu$
		\end{tabular}
		\label{f:predinfo-bg}
		\caption{
		Venn diagram representation \cite[Ch. 2]{CoverThomas} of several information measures for
		stationary random processes. Each circle or oval represents a random
		variable or sequence of random variables relative to time $t=0$. Overlapped areas
		correspond to various mutual informations.
		In (c), the circle represents the `present'. Its total area is
		$H(X_0)=H(1)=\rho_\mu+r_\mu+b_\mu$, where $\rho_\mu$ is the multi-information
		rate, $r_\mu$ is the residual entropy rate, and $b_\mu$ is the predictive
		information rate. The entropy rate is $h_\mu = r_\mu+b_\mu$.
		}
	\end{figure}
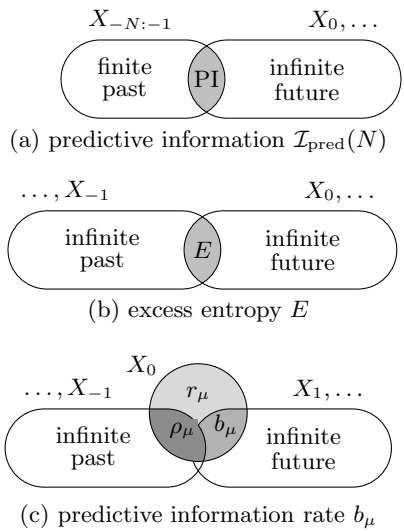

	In \cite{AbdallahPlumbley2009}, we introduced 
	the \emph{predictive information rate}
	(PIR), which is the average information 
	in one observation about the infinite future given the infinite past.
	If $\past{X}_t=(\ldots,X_{t-2},X_{t-1})$ denotes the variables
	before time $t$, 
	and $\fut{X}_t = (X_{t+1},X_{t+2},\ldots)$ denotes
	those after $t$,
	the PIR is defined as a conditional mutual information: 
	\begin{equation}
		\label{eq:PIR}
		\IXZ_t \define I(X_t;\fut{X}_t|\past{X}_t) = H(\fut{X}_t|\past{X}_t) - H(\fut{X}_t|X_t,\past{X}_t).
	\end{equation}
	Equation \eqrf{PIR} can be read as the average reduction
	in uncertainty about the future on learning $X_t$, given the past. 
	Due to the symmetry of the mutual information, it can also be written
	as $\IXZ_t = H(X_t|\past{X}_t) - H(X_t|\fut{X}_t,\past{X}_t)$.
	$H(X_t|\past{X}_t)$ 
	is the entropy rate $h_\mu$, but $H(X_t|\fut{X}_t,\past{X}_t)$ is 
	a quantity that does not appear to be have been considered by other authors yet.
	It is the conditional entropy of one variable given \emph{all} the others 
	in the sequence, future as well as past.
	We call this the \emph{residual entropy rate}
	$r_\mu$, and define it as a limit:
	\begin{equation}
		\label{eq:residual-entropy-rate}
		r_\mu \define \lim_{N\tends\infty} H(X_{-N..N}) - H(X_{-N..-1},X_{1..N}).
	\end{equation}
	The second term, $H(X_{-N..-1},X_{1..N})$, 
	is the joint entropy of two non-adjacent blocks with a 
	gap between them, and cannot be expressed as a function of block entropies alone.
	If we 
	let $b_\mu$ denote the shift-invariant PIR, 
	then
	$b_\mu = h_\mu - r_\mu$ (see \Figrf{predinfo-bg}).


	Many of the measures reviewed above were 
	intended as measures of `complexity', 
	a quality that is somewhat open to interpretation 
	\cite{Bennett1990,FeldmanCrutchfield1998}. 
	It is generally agreed, however, that
	complexity should be low for systems that are deterministic
	or easy to compute or predict---`ordered'---and low
	for systems that a completely random and unpredictable---`disordered'.
%
	The PIR 
	satisfies these conditions 
	without being `over-universal' in the sense of Crutchfield \etal
	\cite{FeldmanCrutchfield1998,CrutchfieldFeldmanShalizi2000}:
	it is not simply a function of entropy or entropy rate
	that fails to distinguish between the different strengths of temporal
	dependency that can be exhibited by systems 
	at a given level of entropy. 
	In our analysis of Markov chains \cite{AbdallahPlumbley2009}, we found
	that processes which maximise the PIR do not maximise the
	multi-information rate $\rho_\mu$ (or the excess entropy,
	which is the same in this case), but do
	have a certain kind of partial predictability
	that requires the observer continually to pay attention to the most
	recent observations
	in order to make optimal predictions. 
	And so, while Crutchfield \etal make a compelling case for the excess
	entropy $E$ and their statistical complexity $C_\mu$ as measures of 
	complexity,
	there is still room to suggest that the PIR captures
	a different and non trivial aspect of temporal dependency structure
	not previously examined.

	\section{Binding information}
	\label{s:binding-info}

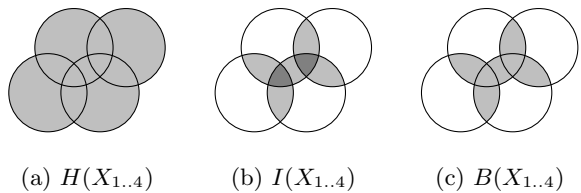
\begin{figure}
		\label{f:multi-binding}
		\centering
		\pgfsetxvec{\pgfpoint{4mm}{0mm}}%
		\pgfsetyvec{\pgfpoint{0mm}{4mm}}%
		\newcommand\colsep{\qquad}
		\newcommand{\rvar}[3]{\coordinate (#3) at (0.866 * #1,#2);}%
		\newcommand{\vars}{\rvar{0}{0}{A}\rvar{1}{1.5}{B}\rvar{2}{0}{C}\rvar{3}{1.5}{D}}%
		\newcommand{\circo}{circle (1.3)}%
		\newcommand{\olap}[2]{%
			\begin{scope}
				\foreach \nd in {#2} {\clip (\nd) \circo;}
				\fill[#1] (-1,-1.2) rectangle (3.5,2.6);
			\end{scope}}
		\def\olapl(#1,#2){\olap{lightgray}{#1,#2}}%
		\def\olapd(#1,#2,#3){\olap{gray}{#1,#2,#3}}%
		\begin{tabular}{c@{\colsep}c@{\colsep}c}
			\begin{tikzpicture}
				\vars%
				\foreach \nd in {A,B,C,D} {\fill[lightgray] (\nd) \circo;}
				\foreach \nd in {A,B,C,D} {\draw (\nd) \circo;}
			\end{tikzpicture} 
			&
			\begin{tikzpicture}
				\vars%
				\olapl(A,B)\olapl(B,C)%
				\olapl(A,C)\olapl(B,D)\olapl(C,D)%
				\olapd(A,B,C)\olapd(B,C,D)%
				\foreach \nd in {A,B,C,D} {\draw (\nd) \circo;}
			\end{tikzpicture} 
			&
			\begin{tikzpicture}
				\vars%
				\olapl(A,B)\olapl(B,C)%
				\olapl(A,C)\olapl(B,D)%
				\olapl(C,D)%
				\foreach \nd in {A,B,C,D} {\draw (\nd) \circo;}
			\end{tikzpicture} 
			\\[1em] 
			\iftwocol{%
				(a) $H(X_{1..4})$ & (b) $I(X_{1..4})$ & (c) $B(X_{1..4})$ 
			}{%
				(a) Joint entropy & (b) Multi-information & (c) Binding information
			}%
		\end{tabular}
		\caption{Illustration of binding information
		as compared with multi-information for a set of four random variables. 
		In each case, the quantity
		is represented by the total amount of black ink, as it were,
		in the shaded parts of the diagram. 
		Whereas the multi-information counts the
		multiply-overlapped areas multiple times, the binding information counts
		each overlapped areas just once.
		}
	\end{figure}

	If the PIR rate is accumulated over
	successive time steps, 
	a quantity which we
	call the \emph{binding information} is obtained. 
	To proceed,
	we first reformulate the infinite sequence
	PIR \eqrf{PIR} so that it becomes applicable to a \emph{finite} sequence of random variables
	$(X_1,\ldots,X_N)$:
	\begin{equation}
		\label{eq:PIR-finite}
		\IXZ_t(X_{1..N}) = I(X_{t};X_{(t+1)..N}|X_{1..(t-1)}),
	\end{equation}
	Note that this is no longer shift-invariant and may depend on $t$.
	The binding information, then, is the sum
	\begin{equation}
		B(X_{1..N}) = \sum_{t\in 1..N} \IXZ_t(X_{1..N}).
		\label{eq:bindinfo-seq}
	\end{equation}
	Expanding this in terms of entropies and conditional entropies and
	cancelling terms yields
	\begin{equation}
		\label{eq:binding-num}
		B(X_{1..N}) = H(X_{1..N}) - \sum_{t \in 1..N} H(X_t|X_{1..N\setminus \{t\}}).
	\end{equation}
	Like the multi-information, 
	it measures dependencies between random
	variables, but in a different way (see \figrf{multi-binding}). 
	Though the binding information was derived by 
	accumulating the PIR sequentially,
	the result is permutation invariant, suggesting that
	the concept might be applicable to
	arbitrary sets of random variables regardless of 
	their topology. Accordingly, we define the binding information
	as follows:

	\begin{definition}
		If $\{X_\alpha| \alpha \in \A\}$ is 
		set of random variables indexed by a countable set $\A$,
		the binding information is 
		\begin{equation}
			\label{eq:binding-def}
			B(X_{\A}) \define H(X_{\A}) - 
				\sum_{\alpha\in\A} H(X_\alpha|X_{\A\setminus \{\alpha\}}).
		\end{equation}
	\end{definition}

%
%

	Since the binding information can be expressed as a sum of 
	(conditional) mutual informations between sets of random variables
	\eqrf{bindinfo-seq}, it is (a) non-negative and (b) invariant to invertible
	pointwise transformations of the variables; 
	that is, if $Y_{\A}$ is a set of random variables such
	that, $\forall \alpha\in\A$,  $Y_\alpha = f_\alpha(X_\alpha)$ 
	for some invertible functions $f_\alpha$,
	then $B(Y_{\A}) = B(X_{\A})$.

	The binding information is zero for sets of
	independent random variables---the case of complete `disorder'---and
	zero when all variables have zero entropy,
	taking known values and representing a certain kind of `order'. 
	However, it is also possible to obtain low binding information
	for \emph{random} systems which are nonetheless very ordered in a certain way. 
	If each variable $X_\alpha$ is some function of $X_{\alpha'}$
	for all $\alpha' \neq \alpha$, then
	the state of the entire system can be read off
	from any one of its component variables. In this case, it is
	easy to show that  
	$B(X_\A) = H(X_\A) = H(X_\alpha)$ 
	for any $\alpha \in \A$, 
	which, as we will see, is relatively low compared with what is possible
	as soon as $N$ becomes appreciably large.
	Thus, binding information is low for both highly `ordered' and highly
	`disordered' systems, but in this case, `highly ordered' does \emph{not}
	simply mean deterministic or known \latin{a priori}: it means 
	the whole is predictable from the smallest of its parts.


\section{Bounds on binding and multi-information}
\label{s:bounds}

\begin{figure}
	\label{f:mult-bind-bounds}
	\centering
	\newcommand\N{6}
	\newcommand\pt[2]{node [shape=circle,inner sep=3pt,anchor=#1] {\textbf{#2}}}
	\pgfsetxvec{\pgfpoint{5.1mm}{0mm}}%
	\pgfsetyvec{\pgfpoint{0mm}{5.1mm}}%
	\begin{tikzpicture}[>=stealth',dot/.style={shape=circle,fill=black,inner sep=1.5pt}]
		\scriptsize	
		\newcommand\icolor{black!50}%
		\newcommand\bcolor{black}%
		\draw (0,0)--(\N,0);	
		\draw (0,0)--(0,\N)--(1mm,\N);	
		\foreach \y in {1,2,...,\N} \draw (1mm,\y)--(0,\y);
		\foreach \x in {1,2,...,\N} \draw (\x,1mm)--(\x,0);
		\path (1,-0.25em) node [anchor=north] {$1$};
		\path (\N,-0.25em) node [anchor=north] {$N$};
		\path (-0.25em,\N) node [anchor=east] {$N$};
		\path (\N/2,-0.7em) node [anchor=north] {$H$ (bits)};
		\path (-1.5em,\N/2) node [anchor=south,text=\icolor,rotate=90] {$I$ (bits)};
		\path (-0.25em,\N/2) node [anchor=south,text=\bcolor,rotate=90] {$B$ (bits)};
		\draw[\icolor,thick] (0,0) node [dot] (known) {}
				-- node [sloped,anchor=north,pos=0.55] {$I<(N\cmin 1)H$} (1,\N-1) node (giant-bit) [dot] {}
				-- node [sloped,anchor=south,pos=0.65] {$I< N\cmin H$} (\N,0) node (indep) [dot] {};
		\draw[\bcolor,thick] (0,0) 
			-- node [pos=0.35,sloped,anchor=north] {$B< H$} (\N-1,\N-1) node (parity) [dot] {} 
			-- node [sloped,anchor=south,pos=0.45] {$B< (N\cmin 1)(N\cmin H)$} (\N,0);
		\pt{known}{east}{a}
		\path (known) \pt{east}{a};
		\path (giant-bit) \pt{south}{b};
		\path (parity) \pt{south}{c};
		\path (indep) \pt{west}{d};
	\end{tikzpicture}
	\hspace{0.25em}
	\begin{tikzpicture}[>=stealth',dot/.style={shape=circle,fill=black,inner sep=1.5pt}]
		\scriptsize	
		\fill[black!15] (0,0) -- (1,\N-1) -- (\N-1,1) -- cycle;
		\draw (0,0)--(\N,0);	
		\draw (0,0)--(0,\N)--(1mm,\N);	
		\foreach \y in {1,2,...,\N} \draw (1mm,\y)--(0,\y);
		\foreach \x in {1,2,...,\N} \draw (\x,1mm)--(\x,0);
		\path (1,-0.25em) node [anchor=north] {$1$};
		\path (-0.25em,1) node [anchor=east] {$1$};
		\path (\N,-0.25em) node [anchor=north] {$N$};
		\path (-0.25em,\N) node [anchor=east] {$N$};
		\path (\N/2,-0.7em) node [anchor=north] {$B$ (bits)};
		\path (-0.5em,\N/2) node [anchor=south,rotate=90] {$I$ (bits)};
		\draw[thick] (0,\N) -- node [sloped,anchor=south] {$I+B<N$} (\N,0);
		\draw[thick] (0,0) node [dot] (known) {}
				-- node [sloped,anchor=north] {$I< (N\cmin 1)B$} (1,\N-1) node (giant-bit) [dot] {};
		\draw[thick] (0,0) 
			-- node [sloped,anchor=south] {$B<(N\cmin 1)I$} (\N-1,1) node (parity) [dot] {};
		\path (giant-bit) \pt{south west}{b};
		\path (parity) \pt{south west}{c};
		\path (0,0) \pt{east}{a} \pt{north}{d};
	\end{tikzpicture}
	\caption{Constraints on multi-information $I(X_{1..N})$ and
		binding information $B(X_{1..N})$ for a system of $N=\N$ binary
		random variables. The labelled points represent identifiable
		distributions over the $2^N$ states that this system can occupy:
		(a) \emph{known state}, the system is deterministically in one configuration;
		(b) \emph{giant bit}, one of the $P^6_\mathcal{B}$ processes;
		(c) \emph{parity}, the parity processes $\parity{6}{0}$ or $\parity{6}{1}$;
		(d) \emph{independent}, the system of independent unbiased random bits. 
	}
\end{figure}
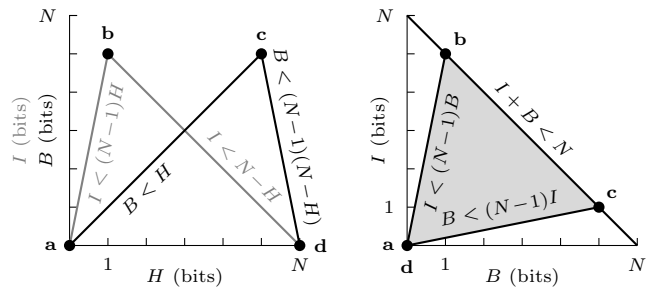

In this section we confine our attention to sets of discrete random
variables taking values in a common alphabet containing $K$ symbols. 
In this case, it is quite straightforward to
derive upper bounds, as functions of the joint entropy, on both the 
multi-information and the binding information, and also
upper bounds on multi-information and binding information as
functions of each other. In \cite{AbdallahPlumbley2010}, we prove the following results:
\begin{theorem}
	\label{th:multi-entropy-bound}
	If $\{X_\alpha| \alpha \in \A\}$ is a set of $N=\card{\A}$ random variables all 
	taking values in a discrete set of cardinality $K$, then the following constraints all hold:
	\begin{align}
		\label{eq:multi-bound-1} I(X_\A) &\leq N \log K - H(X_\A) \\
		\label{eq:multi-bound-2} I(X_\A) &\leq (N - 1)H(X_\A) \\
		B(X_\A) &\leq H(X_\A) \\
		B(X_\A) &\leq (N - 1)(N\log K - H(X_\A)). 
	\end{align}
	Also, $B(X_\A)$ and $I(X_\A)$ are mutually constrained:
	\begin{align}
		\iftwocol{I(X_\A) + B(X_\A) &\leq N \log K.}
		{B(X_\A) + I(X_\A) &\leq N \log K.}
	\end{align}
\end{theorem}
These bounds restrict $I(X_\A)$ and $B(X_\A)$
to two triangular regions of the plane when plotted against the joint entropy
$H(X_\A)$ and are illustrated for $N=6, K=2$ in \figrf{mult-bind-bounds}.
Two more linear bounds were suggested by empirical computations of binding
information and multi-information:
	\begin{align}
		\label{eq:mult-bind-bound1}
		I(X_\A)  &\leq (N-1)B(X_\A) \\
	\text{and}\quad 
		\label{eq:mult-bind-bound2}
		B(X_\A) &\leq (N-1)I(X_\A).
	\end{align}
We have not found a general proof of these inequalities 
for all $N$, but we have constructed a numerical algorithm 
\cite{AbdallahPlumbley2010} that is able to
find proofs for given values of $N$
up to $37$, at which point insufficient numerical precision 
becomes the limiting factor.

\section{Maximising binding information}
\label{s:maximising}

Is
the absolute maximum of $B(X_{1..N}) = (N-1)\log K$ implied by Theorem 
\ref{th:multi-entropy-bound} is attainable, and by what kinds of processes? 
In \cite{AbdallahPlumbley2010} we prove the following:
\begin{theorem}
	If $\{X_1,\ldots,X_N\}$ is a set of discrete random variables
	each taking values in $0..(K\cmin 1)$, then 
	$B(X_{1..N})$ is maximised at $(N\cmin 1)\log_2 K$ bits by the $K$
	`modulo-$K$ processes' $P_{K,m}^N$ for $m \in 0..(K\cmin 1)$,
	under which the probability of a configuration $\vec{x}\in (0..K\cmin 1)^N$ is
	\begin{equation}
		P_{K,m}^N(\vec{x}) = \begin{cases}
			K^{1-N} & \text{if }\modulo{\left( \sum_{i=1}^N x_i \right)}{K} = m,\\
			0 & \text{otherwise.}
		\end{cases}
	\end{equation}
\end{theorem}
When $K=2$ (binary random variables)
the maximal binding information of $N\cmin 1$ bits is
reached by the two `parity' processes:
$P^N_{2,0}$ is the `even' process, which distributes 
uniform probability over all configurations with even parity; 
$P^N_{2,0}$ is the `odd'
process, which distributes uniform probabilities over the complementary set.
The multi-information of the parity processes 
is 1 bit.
%
By contrast, the binary processes which maximise the multi-information
at $N\cmin 1$ bits
are the `giant bit' processes:
the indices $1..N$ are partitioned
into two sets $\mathcal{B}$ and its complement
$\overline{\mathcal{B}} =  1..N\setminus \mathcal{B}$,
and probabilities assigned to configurations $\vec{x}\in\{0,1\}^N$
as follows:
\begin{equation}
	P^N_\mathcal{B}(\vec{x}) = \begin{cases}
		\half &: \text{if } \forall i\in1..N \sps x_i= \mathbb{I}(i\in\mathcal{B}),  \\
		\half &: \text{if } \forall i\in1..N \sps x_i= \mathbb{I}(i\in\overline{\mathcal{B}}), \\
		0     &: \text{otherwise},
	\end{cases}
\end{equation}
where $\mathbb{I}(\cdot)$ is $1$ if 
its argument
is true and $0$ otherwise.
The binding information of these processes is 1 bit.
Thus we see that the processes which maximise the binding information
and the multi-information are quite different in character.

\section{Discussion and conclusions}
\label{s:discussion}

	As noted in \secrf{intro}, Bialek \etal 
	argue that the predictive information $\Ipred(N)$, being
	the sub-extensive component of the entropy, 
	is the unique measure of complexity that satisfies certain
	reasonable desiderata, including transformation invariance
	for continuous-valued variables
	\cite[\S 5.3]{BialekNemenmanTishby2001}. 
	While lack of space precludes a full discussion,
	we note that transformation invariance 
	does \emph{not}, as 
	Bialek \etal state \cite[p. 2450]{BialekNemenmanTishby2001}, 
	demand sub-extensivity:
	binding information \emph{is} transformation
	invariant, since it is a sum of conditional mutual
	informations, and yet it 
	\emph{can} have an extensive component, 
	since its intensive counterpart,  the PIR,
	can have a well-defined value, \eg, in
	stationary Markov chains \cite{AbdallahPlumbley2009}. 



	Measures of statistical dependency are discussed by Studen{\`y} and Vejnarov{\`a},
	\cite[\S 4]{StudenyVejnarova1998},
	who formulate a `level-specific' measure 
	that captures the dependency visible when \emph{fixed size} subsets of
	variables are examined in isolation.
	Studen{\`y} and Vejnarov{\`a} \cite[p. 277]{StudenyVejnarova1998} 
	use the parity process as an 
	example of a random process in which
	the dependence is only visible at the highest level, that is, amongst
	all $N$ variables; if fewer than $N$ variables are examined, they 
	appear to be independent. They note that such processes were
	called `pseudo-independent' by Xiang \etal \cite{XiangWongCercone1996},
	who concluded that standard algorithms for Bayesian network construction
	fail when applied to them. It is intriguing, then, that 
	these are singled out as `most complex'
	according to the binding information criterion.

	To summarise, we have introduced binding information
	as a measure of statistical structure that can be applied to
	any countable set of random variables regardless of
	any topological organisation of the variables.
	Binding information 
	is maximised in finite discrete valued systems
	by the `modulo process'.
	Further results on binding information, 
	and investigations of binding information
	in some specific random processes
	are presented in \cite{AbdallahPlumbley2010}.

	\bibliography{all,c4dm,compsci}
	\comment{
		\begin{acknowledgments}
			This research was supported by EPSRC grant GR/S82213/01.
		\end{acknowledgments}
	}
\end{document}